# The *Z*-Transform of the Levi-Civita Symbol

W. Astar
University of Maryland, Baltimore County(UMBC); Baltimore, Maryland 21250

***Abstract*** - A derivation of the *Z*-transform of the Levi-Civita symbol is carried out up to 3 dimensions, and generalized to higher dimensions. Using its *Z*-transforms, corresponding expressions for the Laplace transform are also found, using Tustin's transform

## 1. Introduction

The Levi-Civita symbol (or epsilon) represents the set {0, 1, -1}, and is generally defined for an *N*-dimensional vector space ($R^N$), as the sign of a permutation of the set of natural numbers {1, 2, 3,···, *N*}. Using the index set or *N*-tuple {$n_1$, $n_2$, $n_3$,···, $n_N$}, in which each index $n_n$ assumes the values of the set {1, 2, 3,···, *N*}, the epsilon is conventionally expressed as

$$\varepsilon_{n_1 n_2 n_3 \cdots n_N} = \begin{cases} +1 & \text{if } \{n_1, n_2, n_3, \cdots, n_N\} \text{ is any even permutation of } \{1, 2, 3, \cdots, N\}, \\ -1 & \text{if } \{n_1, n_2, n_3, \cdots, n_N\} \text{ is any odd permutation of } \{1, 2, 3, \cdots, N\}, \\ 0 & \text{otherwise.} \end{cases} \quad (1.1)$$

In $R^N$, an epsilon generates $N^N$ distinct *N*-tuples, out of which (*N*!) *N*-tuples have entirely dissimilar values, whereas ($N^N$-*N*!) *N*-tuples have at least 2 identical values. Furthermore, the epsilon changes its sign if any 2 values are interchanged in an *N*-tuple [1].

In two-dimensions ($R^2$), the epsilon is defined as

$$\varepsilon_{ij} = \begin{cases} +1 & \text{if } \{i, j\} = \{1, 2\}, \\ -1 & \text{if } \{i, j\} = \{2, 1\}, \\ 0 & \text{if } \{i, j\} \in \{\{1,1\}, \{2,2\}\}. \end{cases} \quad (1.2)$$

In three-dimensions ($R^3$),

$$\varepsilon_{ijk} = \begin{cases} +1 & \text{if } \{i, j, k\} \in \{\{1,2,3\}, \{2,3,1\}, \{3,1,2\}\}, \\ -1 & \text{if } \{i, j, k\} \in \{\{1,3,2\}, \{2,1,3\}, \{3,2,1\}\}, \\ 0 & \text{if } i = j, \text{ or } i = k, \text{ or } j = k. \end{cases} \quad (1.3)$$

For higher dimensional vector spaces, this description of the epsilon becomes increasingly cumbersome, since the number of *N*-tuples for which the epsilon is non-zero increases as *N*!. For instance, continuing with the above approach for the epsilon in four-dimensions ($R^4$) yields

$$\varepsilon_{ijkl} = \begin{cases} +1 & \text{if } \{i,j,k,l\} \in \begin{cases} \{1,2,3,4\},\{1,3,4,2\},\{1,4,2,3\};\{2,1,4,3\},\{2,3,1,4\},\{2,4,3,1\}; \\ \{3,1,2,4\},\{3,2,4,1\},\{3,4,1,2\};\{4,1,3,2\},\{4,2,1,3\},\{4,3,2,1\} \end{cases}, \\ -1 & \text{if } \{i,j,k,l\} \in \begin{cases} \{1,3,2,4\},\{1,4,3,2\},\{1,2,4,3\};\{2,4,1,3\},\{2,1,3,4\},\{2,3,4,1\}; \\ \{3,2,1,4\},\{3,4,2,1\},\{3,1,4,2\};\{4,3,1,2\},\{4,1,2,3\},\{4,2,3,1\} \end{cases}, \\ 0 & \text{otherwise.} \end{cases}$$

(1.4)

In a previous report [2], it was pointed out that this description has more in common with look-up tables, or to sets of conditions, than with an equation that admits values into variables.

There have been attempts at alternative, compact representations of the epsilon, one of which is the following expression, which is in terms of the signum function [1],

$$\varepsilon_{n_1 n_2 n_3 \cdots n_N} = \prod_{N \geq q \geq p \geq 1} \text{sgn}(n_q - n_p); \quad n_1, n_2, n_3, \cdots, n_N \in \{1, 2, 3, \cdots, N\} \quad (1.5)$$

which is an improvement over (1.1), and is especially useful when the number of dimensions exceeds three. However, the signum function is neither an elementary [3, 4], nor an analytical function [5]. Furthermore, (1.5) still incorporates an unresolved condition $N \geq q \geq p \geq 1$ in the form of a multiple inequality.

In a previous report [2], an analytical expression was derived for the epsilon (1.1) in $R^2$, and in $R^3$, and was deduced to be in $R^N$ as

$$\varepsilon_{n_1 n_2 n_3 \cdots n_N} = \prod_{p=1}^{N-1} \prod_{q=1}^{N-p} \frac{n_{N+1-p} - n_q}{N+1-p-q}; \quad n_1, n_2, n_3, \cdots, n_N \in \{1, 2, 3, \cdots, N\}, \quad N > 1. \quad (1.6)$$

This expression was verified in Matlab for $R^5$ [2]. It was also found that the epsilon (1.1) can be re-cast as a generalized, discrete function

$$\varepsilon_{n_1 n_2 n_3 \cdots n_N}(\lambda) = \prod_{p=1}^{N-1} \prod_{q=1}^{N-p} \frac{G(n_{N+1-p}\lambda) - G(n_q \lambda)}{G((N+1-p)\lambda) - G(q\lambda)}; \quad n_1, n_2, n_3, \cdots, n_N \in \{1, 2, 3, \cdots, N\}; N > 1; \lambda \in \mathbf{C}$$

(1.7)

in which $\lambda$ is an arbitrary complex constant. Many elementary and non-elementary functions, such as trigonometric and Bessel functions, were found to satisfy the above the relation, which reproduces (1.1) and (1.6) for the identity function $G(\xi) = \xi$.

Although not so obvious from (1.1), (1.6) should however make the epsilon amenable to various transform techniques, among which is the Z-transform [6 - 8]

$$X(z) = \sum_{n=-\infty}^{\infty} x_n z^{-n} \quad (1.8)$$

which is sometimes termed the bilateral Z-transform. For finite sequences such as the epsilon that are also indexed to *n* being no less than unity, the above expression simplifies to the finite, unilateral Z-transform,

$$X(z) = \sum_{n=-\infty}^{\infty} \left( u_{n-1} - u_{n-N-1} \right) x_n z^{-n} = \sum_{n=1}^{N} x_n z^{-n} \tag{1.9}$$

which effectively renders *x* a causal function. The equation makes use of the discrete, Heaviside step-function, variously defined as [8]

$$u_{n-n_0} = u[n - n_0] = \begin{cases} 0 &, n < n_0 \\ 1 &, n \geq n_0 \end{cases} \tag{1.10}$$

The Z-transform (1.9) is usually resolved using well-known expressions for finite (geometric) series. The result is accompanied by its region of convergence (ROC), which may be derived from the series, and which represents the set of values in the complex plane for which the Z-transform sum is convergent. At a minimum, the ROC is usually expressed as a single inequality in *z*. The ROC is important to specify, because identical Z-transforms can have different ROCs. For the epsilon (1.6), which is always a multi-dimensional quantity, the transform (1.9) may be carried out along a single dimension. However, it should also be possible to carry out the multi-dimensional unilateral Z-transform on (1.6), expressed here as

$$X(z_1, z_2, z_3, \cdots, z_N) = \prod_{i=1}^{N} \sum_{n_i=1}^{N} z_i^{-n_i} [x_{n_1 n_2 n_3 \cdots n_N}], \tag{1.11}$$

with the implication being that the quantity enclosed by the square brackets is the subject of multiple (*N*) single-dimensional operations (1.9) represented by summations over $z_i^{-n_i}$. The epsilon expression (1.6) is not product-wise separable in its dimensions, which can complicate the resultant Z-transform.

Once the Z-transform is found for an *N* dimensional quantity, the Laplace transform can be found using Tustin's bilinear transform of [9]

$$z_n = e^{s_n T} \approx \frac{1 + s_n T/2}{1 - s_n T/2}; \quad n = \{1, 2, 3 \cdots, N\} \tag{1.12}$$

where *T* is conventionally a positive constant, and may be taken as the reciprocal of the sampling frequency used to attain the discrete (digital) function from its continuous counterpart. Thus, $s_n$ has units of frequency or Hertz. The transformation is a method of conformal mapping, and maps the digital, or $z_n$-plane to the analog, or $s_n$-plane. More explicitly, (1.12) maps the perimeter of the unit circle centered at the origin of the $z_n$-plane, to the imaginary axis of the $s_n$-plane. Furthermore, the interior of this unit circle is mapped to the left half of the $s_n$-plane, whereas its exterior, to the right half of the $s_n$-plane, thus preserving stability over the map from the $s_n$-plane to the $z_n$-plane. However, since the mapping (1.12) is non-linear, it results in "frequency warping" [9].

## 2. Two-dimensional vector space

In two-dimensions ($R^2$) or $N = 2$, the Levi-Civita symbol, or the epsilon, is defined as

$$\varepsilon_{n_1 n_2} = \begin{cases} +1 & \text{if } \{n_1, n_2\} = \{1, 2\}, \\ -1 & \text{if } \{n_1, n_2\} = \{2, 1\}, \\ 0 & \text{if } \{n_1, n_2\} \in \{\{1,1\}, \{2,2\}\} \end{cases} \qquad (2.1)$$

but has the analytical expression of (1.6) with $N = 2$:

$$\varepsilon_{n_1 n_2} = n_2 - n_1 \,;\ n_1, n_2 \in \{1, 2\}\,. \qquad (2.2)$$

Its two-dimensional Z-transform (1.11) with $N = 2$, is expressed as

$$E(z_1, z_2) = \sum_{n_2=1}^{2}\sum_{n_1=1}^{2} \varepsilon_{n_1 n_2} z_1^{-n_1} z_2^{-n_2} = \sum_{n_2=1}^{2}\sum_{n_1=1}^{2} (n_2 - n_1) z_1^{-n_1} z_2^{-n_2} \qquad (2.3)$$

which is comprised of a $Z_1$-transform, followed by a $Z_2$-transform. When carried out in this order, the variable $n_2$ is considered a constant under the $Z_1$-transform,

$$E(z_1, z_2) = \sum_{n_2=1}^{2}\left[(n_2 - 1)z_1^{-1} + (n_2 - 2)z_1^{-2}\right]z_2^{-n_2} \qquad (2.4)$$

and is followed by the $Z_2$-transform, which yields the desired Z-transform:

$$E(z_1, z_2) = z_1^{-1}z_2^{-2} - z_1^{-2}z_2^{-1} = \frac{z_1 - z_2}{z_1^2 z_2^2}\,;\ \text{ROC}: \{z_1, z_2\} \in \{\{z_1 > 0\} \times \{z_2 > 0\}\}. \qquad (2.5)$$

It is a two-dimensional transform in $z_1$ and $z_2$, each of which represents a two-dimensional complex plane distinct from the other. The transform exhibits intra-dimensional poles at $z_1 = 0$ and $z_2 = 0$, and an inter-dimensional zero at $z_1 = z_2$, the last of which is analogous to (2.2), which also vanishes for identical indices. Both poles are centered within their respective unit-circles, rendering the transform intra-dimensionally stable. Like the epsilon itself (2.2), its Z-transform (2.5) is not product-wise separable. The ROC is the Cartesian product of the 2 sets in $z_1$ and $z_2$, as shown in (2.5).

The Laplace transform is found by setting $z_n = e^{s_n T}$ in (2.5), or by the bilinear transformation (1.12),

$$E(s_1, s_2) = 4T(s_1 - s_2)\frac{(Ts_1 - 2)(Ts_2 - 2)}{(Ts_1 + 2)^2(Ts_2 + 2)^2}\,. \qquad (2.6)$$

In this case, the intra-dimensional poles have been mapped from zero to $-2/T$ in the left half of the $s_n$-planes, but with their multiplicity intact. However, in addendum to the inter-dimensional zero at $s_1 = s_2$, there now exists intra-dimensional zeros at $2/T$. All the poles and zeros are on their respective real number lines, except the zero at $s_1 = s_2$.

## 3. Three-dimensional vector space

In three-dimensions ($R^3$),

$$\varepsilon_{n_1 n_2 n_3} = \begin{cases} +1 & \text{if } \{n_1, n_2, n_3\} \in \{\{1,2,3\},\{2,3,1\},\{3,1,2\}\}, \\ -1 & \text{if } \{n_1, n_2, n_3\} \in \{\{1,3,2\},\{2,1,3\},\{3,2,1\}\}, \\ 0 & \text{if } n_1 = n_2, \text{ or } n_1 = n_3, \text{ or } n_2 = n_3. \end{cases} \quad (3.1)$$

In [2], the following, alternative analytical expression was derived for (3.1),

$$\varepsilon_{n_1 n_2 n_3} = \frac{1}{2}(n_2 - n_1)(n_3 - n_1)(n_3 - n_2); \quad n_1, n_2, n_3 \in \{1, 2, 3\}. \quad (3.2)$$

It can also be found from (1.6) with $N = 3$. Its three-dimensional Z-transform is found from (1.11) with $N = 3$,

$$E(z_1, z_2, z_3) = \sum_{n_3=1}^{3} \sum_{n_2=1}^{3} \sum_{n_1=1}^{3} \varepsilon_{n_1 n_2 n_3} z_1^{-n_1} z_2^{-n_2} z_3^{-n_3} = \frac{1}{2} \sum_{n_3=1}^{3} \sum_{n_2=1}^{3} \sum_{n_1=1}^{3} (n_2 - n_1)(n_3 - n_1)(n_3 - n_2) z_1^{-n_1} z_2^{-n_2} z_3^{-n_3}$$

(3.3)

Since the variables $n_2$ and $n_3$ are considered as constants under the $Z_1$-transform,

$$E(z_1, z_2, z_3) = \frac{1}{2} \sum_{n_3=1}^{3} z_3^{-n_3} \sum_{n_2=1}^{3} (n_3 - n_2) z_2^{-n_2} \sum_{n_1=1}^{3} \left( n_1^2 - (n_2 + n_3) n_1 + n_2 n_3 \right) z_1^{-n_1} \quad (3.4)$$

whereas $n_1$ is a constant under the $Z_2$-transform, yielding

$$E(z_1, z_2, z_3) = \frac{1}{2} \sum_{n_3=1}^{3} z_3^{-n_3} \left[ \left( n_3 \sum_{n_2=1}^{3} z_2^{-n_2} - \sum_{n_2=1}^{3} n_2 z_2^{-n_2} \right) \sum_{n_1=1}^{3} n_1^2 z_1^{-n_1} + \left( \sum_{n_2=1}^{3} n_2^2 z_2^{-n_2} - n_3^2 \sum_{n_2=1}^{3} z_2^{-n_2} \right) \sum_{n_1=1}^{3} n_1 z_1^{-n_1} \right. \\ \left. + \left( n_3^2 \sum_{n_2=1}^{3} n_2 z_2^{-n_2} - n_3 \sum_{n_2=1}^{3} n_2^2 z_2^{-n_2} \right) \sum_{n_1=1}^{3} z_1^{-n_1} \right].$$

(3.5)

After carrying out the $Z_3$-transform, (3.5) is finally simplified as

$$E(z_1, z_2, z_3) = \frac{1}{2} \left[ \left( \sum_{n_3=1}^{3} n_3 z_3^{-n_3} \sum_{n_2=1}^{3} z_2^{-n_2} - \sum_{n_3=1}^{3} z_3^{-n_3} \sum_{n_2=1}^{3} n_2 z_2^{-n_2} \right) \sum_{n_1=1}^{3} n_1^2 z_1^{-n_1} + \left( \sum_{n_3=1}^{3} z_3^{-n_3} \sum_{n_2=1}^{3} n_2^2 z_2^{-n_2} - \sum_{n_3=1}^{3} n_3^2 z_3^{-n_3} \sum_{n_2=1}^{3} z_2^{-n_2} \right) \sum_{n_1=1}^{3} n_1 z_1^{-n_1} \\ + \left( \sum_{n_3=1}^{3} n_3^2 z_3^{-n_3} \sum_{n_2=1}^{3} n_2 z_2^{-n_2} - \sum_{n_3=1}^{3} n_3 z_3^{-n_3} \sum_{n_2=1}^{3} n_2^2 z_2^{-n_2} \right) \sum_{n_1=1}^{3} z_1^{-n_1} \right].$$

(3.6)

Each summation is of the following form,

$$\sum_{n_q=1}^{3} n_q^p z_q^{-n_q} = \frac{z_q^2 + 2^p z_q + 3^p}{z_q^3} = S(p,q)\,;\; p \in \{0,1,2\}\,,\; q \in \{1,2,3\}\,;\; \text{ROC: } z_q > 0 \;\forall\, q$$

(3.7)

which is a $Z_q$-transform, therefore simplifying (3.6) to

$$E(z_1,z_2,z_3) = \frac{1}{2}\left(\begin{array}{c}\bigl(S(1,3)S(0,2) - S(0,3)S(1,2)\bigr)S(2,1) + \bigl(S(0,3)S(2,2) - S(2,3)S(0,2)\bigr)S(1,1) \\ + \bigl(S(2,3)S(1,2) - S(1,3)S(2,2)\bigr)S(0,1)\end{array}\right).$$

(3.8)

Using the Kronecker delta symbol, the equation can be re-cast as

$$E(z_1,z_2,z_3) = \frac{1}{2}\sum_{m=1}^{3}(-1)^{m-1}\bigl(S(2\delta_{m-1}+2\delta_{m-2}+\delta_{m-3},3)S(\delta_{m-1},2) - S(\delta_{m-1},3)S(2\delta_{m-1}+2\delta_{m-2}+\delta_{m-3},2)\bigr)S(m-1,1).$$

(3.9)

The Kronecker delta symbol is not an elementary function. However, using the following conversion that generalizes various results in [10],

$$\delta_{m-p} = [2\Gamma(m)\cos(p\pi) + m - 2](\Gamma(p)-1) - (p\Gamma(m) - m)\cos(p\pi) + 1\,;\; m,p \in \{1,2,3\}$$

(3.10)

in which $\Gamma(\cdot)$ is the gamma function, then (3.9) can be re-cast in terms of standard functions,

$$E(z_1,z_2,z_3) = \frac{1}{2}\sum_{m=1}^{3}(-1)^{m}\bigl(S(\Gamma(m)-m+1,3)S(3-\Gamma(m),2) - S(3-\Gamma(m),3)S(\Gamma(m)-m+1,2)\bigr)S(m-1,1).$$

(3.11)

After substituting (3.7) for each sum,

$$E(z_1,z_2,z_3) = \sum_{m=1}^{3}\frac{(-1)^m}{2(z_1 z_2 z_3)^3}\left(\begin{array}{c}\bigl(z_3^2 + 2^{\Gamma(m)-m+1}z_3 + 3^{\Gamma(m)-m+1}\bigr)\bigl(z_2^2 + 2^{3-\Gamma(m)}z_2 + 3^{3-\Gamma(m)}\bigr) \\ -\bigl(z_3^2 + 2^{3-\Gamma(m)}z_3 + 3^{3-\Gamma(m)}\bigr)\bigl(z_2^2 + 2^{\Gamma(m)-m+1}z_2 + 3^{\Gamma(m)-m+1}\bigr)\end{array}\right)\bigl(z_1^2 + 2^{m-1}z_1 + 3^{m-1}\bigr)$$

(3.12)

Introducing a 2nd summation in $k$, a more compact expression for the $Z$-transform of (3.1, 2) is finally obtained as

$$E(z_1,z_2,z_3) = \sum_{m=1}^{3}\sum_{k=1}^{2}\frac{(-1)^{m+k}}{2(z_1 z_2 z_3)^3}\left(z_1^2 + \frac{2^m z_1}{2} + \frac{3^m}{3}\right)\left(z_{4-k}^2 + \frac{z_{4-k}}{2^{\Gamma(m)-3}} + \frac{1}{3^{\Gamma(m)-3}}\right)\left(z_{k+1}^2 + \frac{2z_{k+1}}{2^{m-\Gamma(m)}} + \frac{3}{3^{m-\Gamma(m)}}\right);$$

$$\text{ROC}: \{z_1, z_2, z_3\} \in \{\{z_1 > 0\} \times \{z_2 > 0\} \times \{z_3 > 0\}\}$$

(3.13)

which is not product-wise separable in $z_1, z_2,$ and $z_3$.

The Laplace transform of (3.1) or (3.2), can be obtained beginning by setting $z_n = e^{s_n T}$ in (3.7), or by using the bilinear transformation (1.12) on (3.7), formally expressed as

$$R(p,q) = \delta_{z_q - \frac{1+s_q T/2}{1-s_q T/2}} \cdot S(p,q) = \frac{(T_q s_q - 2)\left(2^p (T_q^2 s_q^2 - 4) - 3^p (T_q s_q - 2)^2 - (T_q s_q + 2)^2\right)}{(T_q s_q + 2)^3}$$

(3.14)

which can be simplified to

$$R(p,q) = \sum_{r=1}^{3} \left(\frac{2 - T_q s_q}{2 + T_q s_q}\right)^r r^p ; \quad p \in \{0,1,2\}, \quad q \in \{1,2,3\}.$$

(3.15)

After reducing (3.11) to the compact form

$$E(z_1, z_2, z_3) = \frac{1}{2} \sum_{m=1}^{3} \sum_{k=1}^{2} (-1)^{k+m} S(m-1,1) S(\Gamma(m) - m + 1, k+1) S(3 - \Gamma(m), 4 - k)$$

(3.16)

the Laplace transform of (3.1, 2) is found to be

$$E(s_1, s_2, s_3) = \prod_{q=1}^{3} \delta_{z_q - \frac{1+s_q T/2}{1-s_q T/2}} E(z_1, z_2, z_3) = \frac{1}{2} \sum_{m=1}^{3} \sum_{k=1}^{2} (-1)^{k+m} R(m-1,1) R(\Gamma(m) - m + 1, k+1) R(3 - \Gamma(m), 4 - k).$$

(3.17)

The RHS is actually a function of the Laplace variables $s_n$ due to the use of (3.15). Eq. (3.17) yields, after the substitution of (3.15) for each of the 3 terms in the summand,

$$E(s_1, s_2, s_3) = \frac{1}{2} \sum_{m=1}^{3} \sum_{k=1}^{2} \sum_{r_1=1}^{3} \sum_{r_2=1}^{3} \sum_{r_3=1}^{3} \left(\frac{2 - T_1 s_1}{2 + T_1 s_1}\right)^{r_1} \left(\frac{2 - T_{k+1} s_{k+1}}{2 + T_{k+1} s_{k+1}}\right)^{r_2} \left(\frac{2 - T_{4-k} s_{4-k}}{2 + T_{4-k} s_{4-k}}\right)^{r_3} (-1)^{k+m} r_1^{m-1} r_2^{\Gamma(m)-m+1} r_3^{3-\Gamma(m)}$$

(3.18)

The intra-dimensional poles and the zeros are all on the real number line in each of the three *s*-planes, as was the case for the Laplace transform of the two-dimensional epsilon analyzed in the previous section.

## 4. Summary and conclusions

In $R^2$, the Levi-Civita symbol, or the epsilon, defined as

$$\varepsilon_{n_1 n_2} = n_2 - n_1 = \begin{cases} +1 & \text{if } \{n_1, n_2\} = \{1,2\}, \\ -1 & \text{if } \{n_1, n_2\} = \{2,1\}, \\ 0 & \text{if } \{n_1, n_2\} \in \{\{1,1\},\{2,2\}\} \end{cases} \quad (4.1)$$

was found to have the Z-transform of

$$E(z_1, z_2) = \sum_{n_2=1}^{2}\sum_{n_1=1}^{2}(n_2 - n_1)z_1^{-n_1}z_2^{-n_2} = \frac{z_1 - z_2}{z_1^2 z_2^2}; \quad \text{ROC}: \{z_1, z_2\} \in \{\{z_1 > 0\} \times \{z_2 > 0\}\} \quad (4.2)$$

However, the Z-transform in this case can also be expressed as the determinant of a 2 x 2 matrix using the alternative resolution of

$$E(z_1, z_2) = \sum_{n_1=1}^{2} z_1^{-n_1} \sum_{n_2=1}^{2} n_2 z_2^{-n_2} - \sum_{n_1=1}^{2} n_1 z_1^{-n_1} \sum_{n_2=1}^{2} z_2^{-n_2} = \begin{vmatrix} S_2(0,1) & S_2(0,2) \\ S_2(1,1) & S_2(1,2) \end{vmatrix} \quad (4.3)$$

which is obtained with the observation that

$$\sum_{n_q=1}^{2} n_q^p z_q^{-n_q} = z_q^{-1} + 2^p z_q^{-2} = \frac{z_q + 2^p}{z_q^2} = S_2(p,q); \quad p \in \{0,1\}, \ q \in \{1,2\}; \quad \text{ROC}: z_q > 0 \ \forall \ q \quad (4.4)$$

In $R^3$, the epsilon is defined as

$$\varepsilon_{n_1 n_2 n_3} = \frac{1}{2}(n_2 - n_1)(n_3 - n_1)(n_3 - n_2) = \begin{cases} +1 & \text{if } \{n_1, n_2, n_3\} \in \{\{1,2,3\},\{2,3,1\},\{3,1,2\}\}, \\ -1 & \text{if } \{n_1, n_2, n_3\} \in \{\{1,3,2\},\{2,1,3\},\{3,2,1\}\}, \\ 0 & \text{if } n_1 = n_2, \text{ or } n_1 = n_3, \text{ or } n_2 = n_3. \end{cases} \quad (4.5)$$

which has the Z-transform

$$E(z_1, z_2, z_3) = \sum_{m=1}^{3}\sum_{k=1}^{2} \frac{(-1)^{m+k}}{2(z_1 z_2 z_3)^3}\left(z_1^2 + \frac{2^m z_1}{2} + \frac{3^m}{3}\right)\left(z_{4-k}^2 + \frac{z_{4-k}}{2^{\Gamma(m)-3}} + \frac{1}{3^{\Gamma(m)-3}}\right)\left(z_{k+1}^2 + \frac{2z_{k+1}}{2^{m-\Gamma(m)}} + \frac{3}{3^{m-\Gamma(m)}}\right);$$

$$\text{ROC}: \{z_1, z_2, z_3\} \in \{\{z_1 > 0\} \times \{z_2 > 0\} \times \{z_3 > 0\}\} \quad (4.6)$$

but is also expressible as the scaled determinant of a 3 x 3 matrix

$$E(z_1, z_2, z_3) = \frac{1}{2} \begin{pmatrix} (S_3(2,3)S_3(1,2) - S_3(1,3)S_3(2,2))S_3(0,1) \\ + (S_3(0,3)S_3(2,2) - S_3(2,3)S_3(0,2))S_3(1,1) \\ + (S_3(1,3)S_3(0,2) - S_3(0,3)S_3(1,2))S_3(2,1) \end{pmatrix} = \frac{1}{2} \begin{vmatrix} S_3(0,1) & S_3(0,2) & S_3(0,3) \\ S_3(1,1) & S_3(1,2) & S_3(1,3) \\ S_3(2,1) & S_3(2,2) & S_3(2,3) \end{vmatrix}$$

(4.7)

where

$$S_3(p,q) = \sum_{n_q=1}^{3} n_q^p z_q^{-n_q} = \frac{z_q^2 + 2^p z_q + 3^p}{z_q^3}; \quad p \in \{0,1,2\}, \quad q \in \{1,2,3\}; \quad \text{ROC: } z_q > 0 \quad \forall q$$

(4.8)

in which the sum $S$ is subscripted with a 3 to distinguish it from that (4.4) used for the epsilon in $R^2$.

Based on the above observations, it is concluded that the Z-transform for the epsilon in $R^N$ (1.6), re-expressed here as

$$\varepsilon_{n_1 n_2 n_3 \cdots n_N} = \frac{\prod_{p=1}^{N-1} \prod_{q=1}^{N-p} (n_{N+1-p} - n_q)}{\prod_{p=1}^{N-1} \prod_{q=1}^{N-p} (N+1-p-q)}; \quad n_1, n_2, n_3, \cdots, n_N \in \{1,2,3,\cdots,N\}, \quad N > 1,$$

(4.9)

and given by (1.11)

$$E(z_1, z_2, \cdots, z_N) = \prod_{i=1}^{N} \sum_{n_i=1}^{N} z_i^{-n_i} [\varepsilon_{n_1 n_2 n_3 \cdots n_N}],$$

(4.10)

can be expressed as the scaled determinant of a $N \times N$ matrix

$$E(z_1, z_2, \cdots, z_N) = \frac{1}{\prod_{p=1}^{N-1} \prod_{q=1}^{N-p} (N+1-p-q)} \begin{vmatrix} S_N(0,1) & S_N(0,2) & S_N(0,3) & \cdots & S_N(0,N-1) \\ S_N(1,1) & S_N(1,2) & S_N(1,3) & \cdots & S_N(1,N-1) \\ S_N(2,1) & S_N(2,2) & S_N(2,3) & \cdots & S_N(2,N-1) \\ \vdots & \vdots & \vdots & \ddots & \vdots \\ S_N(N-1,1) & S_N(N-1,2) & S_N(N-1,3) & \cdots & S_N(N-1,N-1) \end{vmatrix}$$

(4.11)

where $S_N(p,q)$ is the $Z_q$-transform of $n_q^p \left[ u_{n_q-1} - u_{n_q-N-1} \right]$, and is given by

$$S_N(p,q) = \sum_{n_q=1}^{N} n_q^p z_q^{-n_q} = \frac{\sum_{r=1}^{N} r^p z_q^{N-r}}{z_q^N}; \quad p \in \{0,1,2,\cdots,N-1\}, \quad q \in \{1,2,\cdots,N\}; \quad \text{ROC: } z_q > 0 \quad \forall q.$$

(4.12)

In $R^2$, the Laplace transform of the epsilon (4.1) was found to be

$$E(s_1, s_2) = 4T(s_1 - s_2)\frac{(Ts_1 - 2)(Ts_2 - 2)}{(Ts_1 + 2)^2 (Ts_2 + 2)^2}. \tag{4.13}$$

It can also be obtained formally from (4.3, 4) by the operation

$$E(s_1, s_2) = \delta_{z_q - \frac{1+s_q T/2}{1-s_q T/2}} \cdot E(z_1, z_2) = \begin{vmatrix} R_2(0,1) & R_2(0,2) \\ R_2(1,1) & R_2(1,2) \end{vmatrix} \tag{4.14}$$

where the converted sum

$$R_2(p,q) = \delta_{z_q - \frac{1+s_q T/2}{1-s_q T/2}} \cdot S_2(p,q) = \sum_{r=1}^{2} \left(\frac{2 - T_q s_q}{2 + T_q s_q}\right)^r r^p \; ; \; p \in \{0,1\}, \; q \in \{1,2\}. \tag{4.15}$$

The Laplace transform of the epsilon in $R^3$ (4.5) may be expressed using (4.7, 8), as

$$E(s_1, s_2, s_3) = \delta_{z_q - \frac{1+s_q T/2}{1-s_q T/2}} \cdot E(z_1, z_2, z_3) = \frac{1}{2}\begin{vmatrix} R_3(0,1) & R_3(0,2) & R_3(0,3) \\ R_3(1,1) & R_3(1,2) & R_3(1,3) \\ R_3(2,1) & R_3(2,2) & R_3(2,3) \end{vmatrix} \tag{4.16}$$

where

$$R_3(p,q) = \delta_{z_q - \frac{1+s_q T/2}{1-s_q T/2}} \cdot S_3(p,q) = \sum_{r=1}^{3} \left(\frac{2 - T_q s_q}{2 + T_q s_q}\right)^r r^p \; ; \; p \in \{0,1,2\}, \; q \in \{1,2,3\}. \tag{4.17}$$

In general then, it can be concluded that for the epsilon in $R^N$ (1.6),

$$E(s_1, s_2, \cdots, s_N) = \delta_{z_q - \frac{1+s_q T/2}{1-s_q T/2}} \cdot E(z_1, z_2, \cdots, z_N) = \frac{\begin{vmatrix} R_N(0,1) & R_N(0,2) & R_N(0,3) & \cdots & R_N(0,N-1) \\ R_N(1,1) & R_N(1,2) & R_N(1,3) & \cdots & R_N(1,N-1) \\ R_N(2,1) & R_N(2,2) & R_N(2,3) & \cdots & R_N(2,N-1) \\ \vdots & \vdots & \vdots & \ddots & \vdots \\ R_N(N-1,1) & R_N(N-1,2) & R_N(N-1,3) & \cdots & R_N(N-1,N-1) \end{vmatrix}}{\prod_{p=1}^{N-1}\prod_{q=1}^{N-p}(N+1-p-q)}$$

(4.18)

where $R_N$ is the Laplace transform of $n_q^p \left[ u_{n_q - 1} - u_{n_q - N - 1} \right]$,

$$R_N(p,q) = \delta_{z_q - \frac{1+s_q T/2}{1-s_q T/2}} \cdot S_N(p,q) = \sum_{r=1}^{N} \left(\frac{2 - T_q s_q}{2 + T_q s_q}\right)^r r^p \; ; \; p \in \{0,1,2,\cdots,N-1\}, \; q \in \{1,2,\cdots,N\}. \tag{4.19}$$

## Acknowledgements

www.WolframAlpha.com was helpful in verifying some of the expressions derived for the Z-transforms and the Laplace transforms. This document was generated in Microsoft Word, and converted prior to the upload using www.freepdfconvert.com